\documentclass{article}

\usepackage[english]{babel}
\usepackage{amsmath,amssymb,amsthm,graphicx}
\usepackage[utf8]{inputenc}

\newtheorem{theorem}{Theorem}

\newtheorem{corollary}[theorem]{Corollary}
\theoremstyle{plain}

\newtheorem{observation}[theorem]{Observation}

\newtheorem{conjecture}[theorem]{Conjecture}

\usepackage{enumerate}
\usepackage{subcaption}
\usepackage{hyperref}
\usepackage{url}

\usepackage[
    natbib = true,
    backend=biber,
    isbn=false,
    url=false,
    doi=false,
    maxbibnames=5,
    minbibnames=5,
    eprint=false,
    style=numeric,
    sorting=nyt,
    sortcites = true
]{biblatex}
\bibliography{xnum}
\DeclareNameAlias{author}{last-first}

\usepackage{xcolor}

\makeatletter

\providecommand{\doi}[1]{%
\begingroup
\let\bibinfo\@secondoftwo
\urlstyle{rm}%
\href{http://dx.doi.org/#1}{%
doi:\discretionary{}{}{}%
\nolinkurl{#1}%
}%
\endgroup
}

\makeatother

\usepackage{etoolbox}
\apptocmd{\sloppy}{\hbadness 10000\relax}{}{}

\newcommand{\crn}{\operatorname{cr}}
\newcommand{\ignore}[1]{}

\title{Graphs with at most one crossing}

\usepackage{authblk}

\author[1]{André C. Silva\thanks{andre.silva@ic.unicamp.br}}
\author[2]{Alan Arroyo\thanks{alanmarcelo.arroyoguevara@ist.ac.at}}
\author[3]{R. Bruce Richter\thanks{rbruce@uwaterloo.ca}}
\author[1]{Orlando Lee\thanks{lee@ic.uniacamp.br}}

\affil[1]{Instituto de Computação, Universidade Estadual de Campinas, Campinas - SP, 13083-852, Brazil}
\affil[2]{IST Austria, 3400 Klosterneuburg, Austria}
\affil[3]{Department of Combinatorics \& Optimization, University of Waterloo, Ont. Waterloo N2L 3G1, Canada}

\begin{document}

\maketitle

\begin{abstract}
	The crossing number of a graph $G$ is the least number of crossings over all possible drawings of $G$. We present a structural characterization of graphs with crossing number one.
\end{abstract}

\newcommand{\janFifteen}[1]{\textcolor{red}{#1}}
\newcommand{\janTwenty}[1]{\textcolor{blue}{#1}}
\newcommand{\janTwentythree}[1]{\textcolor{magenta}{#1}}
\newcommand{\apriltwoone}[1]{\textcolor{red}{#1}}
\newcommand{\april}[1]{\textcolor{magenta}{\Large #1}}

\section{Introduction}

For a graph $G$, we let $V(G)$ and $E(G)$ denote its vertex set and edge set, respectively. Our graphs may have multiple edges but, for simplicity, we assume they have no loops.

We assume the reader is familiar with the concept of drawings of graphs in the plane or sphere. We make no distinction between the elements of a graph and their representations in the drawing. Let $D$ be a drawing of a graph $G$.  If $H$ is a subgraph of $G$, we use the notation $D[H]$ to denote the drawing of $H$ obtained from $D$ by deleting the corresponding vertices and edges that are not in $H$. 

Let $\crn(D)$ denote the number of crossings of $D$. The \textit{crossing number} $\crn(G)$ of $G$ is the least number of crossings over all possible drawings of $G$. A drawing $D$ of $G$ is \textit{optimal} if $\crn(D)=\crn(G)$. We note that a graph $G$ is planar if and only if $\crn(G)=0$. An optimal drawing of a planar graph is an \textit{embedding}.

In the context of studying crossing numbers, we can interpret Kuratowski's classic characterization of planar graphs as: a graph has crossing number at least one if and only if it contains a subdivision of  $K_5$ or $K_{3,3}$ (we shall refer to the subdivisions of $K_{3,3}$ and $K_5$  as \textit{Kuratowski graphs}). Following the same spirit, we answer the question of when does a graph have crossing number at least 2?

We answer with Theorem \ref{thm:1cgc}, a characterization of the crossing pairs of a graph. A pair $\{e,f\}$ of edges in a graph $G$ is a {\em crossing pair\/} of $G$ if there exists a drawing $D$ of $G$ with $\crn(D)=1$ (we refer to $D$ as a \textit{$1$-drawing} of $G$) in which $e$ and $f$ cross.  Clearly, for a nonplanar $G$, $\crn(G)=1$ if and only if $G$ has a crossing pair and $\crn(G) \geq 2$ otherwise. This characterization is an extension of a result of Arroyo and Richter ~\cite{cgcr2} (Theorem~\ref{thm:p4cr2}, below). 

In Section~\ref{sec:cr1rl}, we present some related work. Section~\ref{sec:ocp} expands on some properties of crossing pairs and details our characterization in Theorem~\ref{thm:1cgc}. In Section~\ref{sec:prel}, we describe the results and notation used for proving  Theorem~\ref{thm:1cgc}; the  proof  is in  Section~\ref{sec:p1cgc}. The last section contains some remarks on Theorem~\ref{thm:1cgc}.

\section{Related work}
\label{sec:cr1rl}

The problem of characterizing graphs with crossing number at least two was already studied by Arroyo and Richter~\cite{cgcr2} in the context of \textit{peripherally 4-connected graphs}.

A graph $G$ is \textit{peripherally 4-connected} if $G$ is 3-connected and, for every vertex 3-cut $X$ of $G$, and, for any partition of the components of $G-X$ into two non-null subgraphs $H$ and $K$, at least one of $H$ or $K$ has just one vertex.
Two edges $e=x_1y_1$ and $f=x_2y_2$ are \textit{linked} if either $e,f$ are incident with a common vertex or there is a 3-cut $X$ in $G$ such that $X \subset \{x_1,y_1,x_2,y_2\}$ and the vertex in $\{x_1,y_1,x_2,y_2\} \backslash X$ induces a trivial component of $G-X$. Otherwise, $e$ and $f$ are \textit{unlinked}. Two edges $e$ and $f$ of $G$ are \textit{separated by cycles} if there exists two (vertex-)disjoint cycles $C_e$ and $C_f$ in $G$ with $e \in E(C_e)$ and $f \in E(C_f)$. 

\begin{theorem}\label{thm:p4cr2}
	\cite{cgcr2} A peripherally 4-connected nonplanar graph $G$ has crossing number at least two if and only if every pair of unlinked edges $\{e,f\}$ in $G$ is separated by cycles.
\end{theorem} 

We note that if we drop the connectivity requirement, then the converse of Theorem~\ref{thm:p4cr2} is no longer true. Consider $G = K_{3,4}$. It is 3-connected but the cut consisting of the part of size 3 shows that it is not peripherally 4-connected. We know that $\crn(K_{3,4}) =2$~\cite{zarankiewicz}; however, no pair of disjoint edges is separated by cycles, as we need at least 8 vertices for 2 disjoint cycles.

There is also some work on a characterization of line graphs that have crossing number one. The line graph $L(G)$ of a graph $G$ is a graph with vertex set $E(G)$ and $a,b \in E(G)$ are adjacent in $L(G)$ if and only if $a,b$ share a common vertex in $G$. Let $\Delta(G)$ denote the maximum degree of a graph $G$. 

Since the edges incident with a vertex of degree $d$ in $G$ induce a complete graph $K_d$ in $L(G)$, a vertex of degree $6$ in $G$ implies $\crn(L(G))\geq \crn(K_6)=3$. Similarly, if $G$ has two vertices with degree $5$, then $\crn(L(G))\geq 2$. Therefore, if $\crn(L(G))\leq 1$, then $G$ has at most one vertex of degree $5$ and $\Delta(L(G)) < 8$.

Kulli, Akka and Beineke~\cite{lgcn1} characterized planar graphs whose line graph has crossing number one, while Jendrol' and Klevu{s}\v{c}~\cite{olgcn1} completed the characterization by including nonplanar graphs:

\begin{theorem} 
	\cite{lgcn1}  For every planar graph $G$, we have $\crn(L(G))=1$ if and only if either:
	\begin{enumerate}[(1)]
		\item $\Delta(G) = 4$ and there is a unique non-cut-vertex of degree 4, or
		\item $\Delta(G)=5$, every vertex of degree 4 is a cut vertex, and there is a unique vertex of degree 5 with at most 3 edges in any block. 
	\end{enumerate}
\end{theorem}

\begin{theorem}
	\cite{olgcn1} For a nonplanar graph $G$, we have $\crn(L(G))=1$ if and only if the following conditions hold:
	\begin{enumerate}[(1)]
		\item $\crn(G)=1$,
		\item $\Delta(G) \leq 4$, and every vertex of degree 4 is a cut vertex of $G$, and 
		\item there exists a drawing of $G$ in the plane with exactly one crossing in which each crossed edge is incident with a vertex of degree 2. 
	\end{enumerate}
\end{theorem}

Akka, Jendrol, Kle\v{s}\v{c}, and Panshetty~\cite{olgcn2} obtained a characterization of planar graphs whose line graph has crossing number two.

A graph $G$ is \textit{k-crossing-critical} if $\crn(G) \geq k$ and every proper subgraph $H$ of $G$ has $\crn(H) < k$. The 1-crossing-critical graphs are exactly the Kuratowski graphs. We note that a graph with crossing number at least 2 contains a 2-crossing-critical graph as a subgraph.

A great deal of attention has been given to 2-crossing-critical graphs~\cite{ocnls,lnpga,cccg,cgcn2,cn,ifccg,bigpaper}.  For a positive integer $n \geq 3$, the \textit{Möbius Ladder} $V_{2n}$ on $2n$ vertices is the graph obtained from a $2n$-cycle by joining vertices with distance $n$ in the cycle. Bokal, Oporowski, Richter and Salazar~\cite{bigpaper} characterized all 2-crossing-critical graphs that:
are not 3-connected;  are 3-connected and have $V_{10}$ as a minor; or are 3-connected and do not have $V_8$ as a minor.
They also showed  that there exists only finitely many 3-connected 2-crossing-critical graphs with no $V_{10}$ minor. 

It remains to characterize or enumerate all the 3-connected 2-crossing-critical graphs that have $V_8$ but not $V_{10}$ as a minor. We hope this work can help determine these remaining $2$-crossing-critical graphs.

\section{Crossing pairs and statement of the main result}
\label{sec:ocp}

The main point of this section is to introduce our main result Theorem \ref{thm:1cgc}. The lead-up to its statement is an analysis of crossing pairs. 

\begin{figure}
	\centering
	\subcaptionbox{\label{fig:v8b}}
	[.31\linewidth]{\includegraphics[scale=0.95]{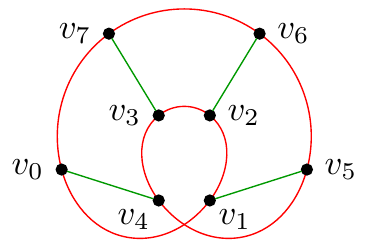}}
	\subcaptionbox{\label{fig:v8b-e}}
	[.31\linewidth]{\includegraphics[scale=0.95]{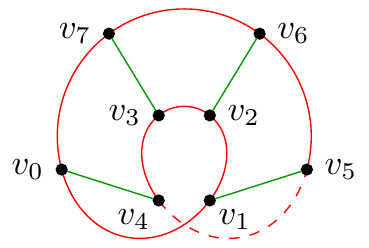}}
	\subcaptionbox{\label{fig:v8bk33}}
	[.31\linewidth]{\includegraphics[scale=0.95]{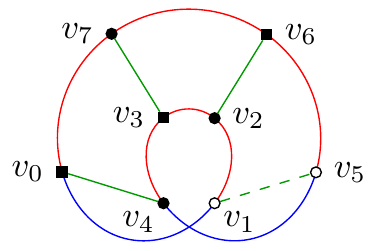}}
		\caption{Removing $v_4v_5$ will result in a planar graph, however removing the edge $v_1v_5$ will result in a subdivision of $K_{3,3}$. The squares and disks represents the parts of the subdivision of $K_{3,3}$.}
		\label{fig:v8}
\end{figure}

The graph in Figure \ref{fig:v8b} is $V_8$. The figure shows a 1-drawing of the $V_8$ where $v_0v_1$ and $v_4v_5$ are crossing. Removing either $v_0v_1$ or $v_4v_5$ results in a planar graph (Figure \ref{fig:v8b-e}). In contrast, as shown in Figure \ref{fig:v8bk33}, removing $v_1v_5$ results in a subdivision of $K_{3,3}$. Since any drawing of $V_8-v_1v_5$ contains a crossing, $v_1v_5$ is not in a crossing pair.

There are three obvious facts about a crossing pair $\{e,f\}$ in a graph $G$. First, as illustrated above,  if $D$ is a $1$-drawing of $G$ in which $e$ and $f$ cross, then $D[G-e]$ and $D[G-f]$ are planar embeddings of $G-e$ and $G-f$, respectively.

Second, if $H$ is a Kuratowski subgraph of $G$, then $e$ and $f$ are both in $H$ and make the unique crossing of $D[H]$. That is, $\{e,f\}$ is a crossing pair in $H$. This is readily seen to be equivalent to the assertion that $e$ and $f$ are not in either the same branch or adjacent branches of $H$. (A \textit{branch} in $H$ is a path $P$ joining two vertices with degrees different from $2$ such that all internal vertices of $P$ have degree $2$ in $H$.)

\begin{figure}
	\centering
	\subcaptionbox{\label{fig:nsepex}}
	[.31\linewidth]{\includegraphics[scale=1]{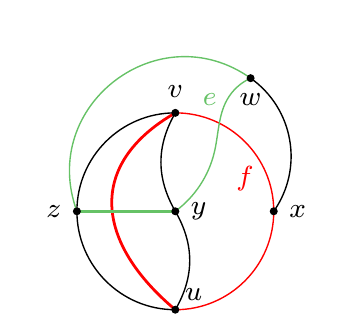}}
	\subcaptionbox{\label{fig:nsepex2}}
	[.31\linewidth]{\includegraphics[scale=1]{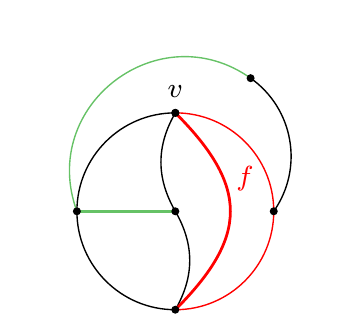}}
	\subcaptionbox{\label{fig:nsepex3}}
	[.31\linewidth]{\includegraphics[scale=1]{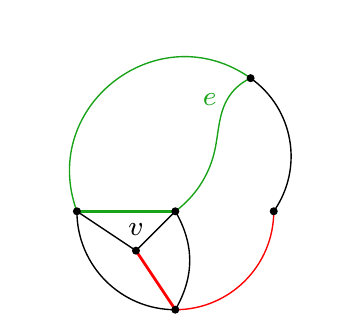}}
		\caption{The pair of edges $e,f$ is separated by the highlighted cycles, however removing either will result in a planar graph.}
		\label{fig:nsepexa}
\end{figure}

Third, since vertex-disjoint cycles cross an even number of times in a drawing, $e$ and $f$ are not separated by cycles in $G$ (separated by cycles is defined just before Theorem \ref{thm:p4cr2}). Figure \ref{fig:nsepexa} shows edges $e$ and $f$ that are separated by cycles such that $G-e$ and $G-f$ are both planar. 

The second fact is easily seen to imply the first: If, for every Kuratowski subgraph $H$ of $G$, $\{e,f\}$ is a crossing pair of $H$, then $G-e$ and $G-f$ are both planar.  However, the first fact does not imply the second: if $e$ and $f$ are in branches of a  Kuratowski graph $H$ sharing a vertex, both $H-e$ and $H-f$ are planar,  but $\{e,f\}$ is not a crossing pair of $H$. In Theorem 4 we describe conditions under which the first and second facts are equivalent.

We are now ready to state our main result, which characterizes crossing pairs in nonplanar graphs in terms of these three facts.

\begin{theorem}\label{thm:1cgc}
Let $G$ be a nonplanar graph and let $e$, $f\in E(G)$. The following are equivalent:
	\begin{enumerate}[\rm(i)]
		\item $\{e,f\}$ is a crossing pair of $G$;
		\item\label{it:two} $\{e,f\}$ is a crossing pair of every Kuratowski subgraph of $G$ and $e$, $f$ are not separated by cycles in $G$; and 
		\item $e$, $f$ are not separated by cycles in $G$, both $G-e$ and $G-f$ are planar and there exists a Kuratowski subgraph $H$ of $G$ such that $\{e,f\}$ is a crossing pair in $H$.
	\end{enumerate}
\end{theorem}

Theorem \ref{thm:1cgc} gives an answer to our original question: what makes a graph have crossing number at least $2$? Let $G$ be a non-planar graph and let $H$ be a Kuratowski subgraph of $G$. Then $\crn(G)\geq 2$ if and only if, for every crossing pair $\{e,f\}$ of $H$, $e$, $f$ are separated by cycles in $G$ or at least one of $G-e$ and $G-f$ is not planar. 

The following example from \v Sir\'a\v n \cite{siran} shows that the two conditions in Theorem \ref{thm:1cgc} (\ref{it:two}) are independent.  Let $(\{u,v,w\},\{x,y,z\})$ be a bipartition of $K_{3,3}$ and let $G$ be the graph $K_{3,3}+\{uv,yz\}$ (see Figure \ref{fig:nsepex}).  As \v Sir\'a\v n notes, there is only one Kuratowski graph contained in $G$.   The edges $ux$ and $wz$ are separated by cycles in $G$.  However, $\{ux,wz\}$ is a crossing pair of every Kuratowski subgraph of $G$.  On the other hand, $\{uy,wz\}$ is also a crossing pair of every Kuratowski subgraph of $G$, but $uy$ and $wz$ are not separated by cycles in $G$.

The proof of Theorem \ref{thm:1cgc} is in Section \ref{sec:p1cgc}. This is preceded by some preparatory work in the next section.

\section{Preliminaries}
\label{sec:prel}

In this section, we present some definitions and results used in the proof of Theorem \ref{thm:1cgc}.

Given a subgraph $H$ of $G$, a path $P$ in $G$ is \textit{$H$-avoiding} or \textit{avoids} $H$ if no edge or internal vertex of $P$ is in $H$. For $u,v \in V(G)$, a $uv$-path in $G$ is a path whose ends are $u$ and $v$.

Let $G$ be a graph and let $H$ be a subgraph of $G$. An \textit{$H$-bridge} $B$ of $G$ is a subgraph of $G$ consisting of either a single edge of $E(G) \setminus E(H)$ with both ends in $H$, or a component $F$ of $G - V(H)$ together with the edges of $G$ with one end in $F$ and another in $H$.
The vertices of $Att(B)=V(H)\cap V(B)$ are the \textit{attachments} of $B$, and the \textit{nucleus} $Nuc(B)$ is $B-Att(B)$. 
Although, the definitions of attachment and nucleus depend on $G$ and on the  subgraph $H$, we omit them, since these are always clear from the context.

Let $C$ be a cycle of $G$. Two distinct $C$-bridges $B_1$ and $B_2$ \textit{overlap} if they have exactly three attachments in common or if there exist distinct vertices $a$, $x$, $b$, $y$ occurring in this cyclic order in $C$ such that $a,b \in Att(B_1)$ and $x,y \in Att(B_2)$. The following simple observations are useful to us.

\begin{observation}\label{obs:overlapping}
Let $C$ be a cycle in a graph $G$ with overlapping $C$-bridges $B_1$ and $B_2$.  
\begin{enumerate}[\rm(i)]\item If $D$ is a planar embedding of $G$, then $B_1$ and $B_2$ are embedded in distinct faces of $D[C]$. 
	\item If $u_1$ and $u_2$ are in the nuclei of $B_1$ and $B_2$, respectively, then $(C\cup B_1\cup B_2)+ u_1u_2$ is not planar. \hfill $\Box$
\end{enumerate}
\end{observation}

For vertices $x$ and $y$ of a graph $G$,  a cycle $C \subseteq G$ \textit{detaches} $x$ from $y$ if there exists two overlapping $C$-bridges, each containing exactly one of $x$ and $y$ in its nucleus.

The vertices $x$ and $y$  are \textit{cofacial} in an embedding $D$ if $x$ and $y$  are incident with a common face of $D$. The following theorem by Tutte (and its slight modification in Corollary~\ref{cor:tuttefc}) is an important tool in the proof of Theorem~\ref{thm:1cgc}.

\begin{theorem}\label{thm:tuttefc}
	(Tutte~\cite{svbc}) Let $G$ be a planar graph and let $x,y\in V(G)$. Then $G$ has an embedding such that $x$ and $y$ are cofacial unless $G$ contains a cycle $C$ that detaches $x$ from $y$.
\hfill\qed\end{theorem}

We need a slight modification of Tutte's result that holds for a vertex $x$ and an edge $f$. A cycle $C$ {\em detaches $x$ from $f$} if there exist overlapping $C$-bridges with one containing $x$ and its nucleus and the other containing $f$. A vertex $x$ 
and an edge $f$ are \textit{cofacial} in an embedding $D$ if $x$ and $f$  are incident with a common face of $D$. 

\begin{corollary}\label{cor:tuttefc}
	Let $G$ be a planar graph, let $x\in V(G)$ and let $f\in E(G)$. Then $G$ has an embedding such that $x$ and $f$ are incident with a common face unless $G$ contains a cycle $C$ that detaches $x$ from $f$. \qed
\end{corollary}

\section{Proof of Theorem \ref{thm:1cgc}}
\label{sec:p1cgc}

In this section we give the proof of Theorem \ref{thm:1cgc}.  
The proof that (i) implies (ii) follows from the second and third facts on the discussion preceding the statement of Theorem \ref{thm:1cgc}.   To see that (ii) implies (iii), by (ii), $G-e$ and $G-f$ both have no Kuratowski subgraphs, so both are planar.  Moreover, $G$ is not planar, so $G$ has a Kuratowski subgraph $H$.  By (ii), $\{e,f\}$ is a crossing pair in $H$ and also $e$ and $f$ are not separated by cycles in $G$.
The rest of this section is devoted to (iii) implies (i).

For a drawing $D$ of a graph $K$, a {\em side} is the closure of a face of $D$.

\begin{proof}  The hypotheses of (iii) are that $e$ and $f$ are not separated by cycles in $G$, $G-e$ and $G-f$ are planar, and there is a Kuratowski subgraph $H$ in which $\{e,f\}$ is a crossing pair. Let $u$ and $v$ be the ends of $e$.

Let $D_H$ be a 1-drawing of $H$ in which $e$ and $f$ cross.
The embedding $D_H[H-e]$ shows that $H-e$ is planar.
Since $H$ is not planar, no embedding of $H-e$ can have both $u$ and $v$ incident with the same face. By Theorem~\ref{thm:tuttefc}, there exists a cycle $C$ of $H-e$, together with distinct overlapping $C$-bridges $B^H_u$ and $B^H_v$ containing $u$ and $v$ in their nuclei, respectively.

Since $B_u^H$ and $B_v^H$ overlap, they must be embedded  in distinct sides of $D_H[C]$ in $D_H[H-e]$. %
Therefore $e$ crosses at least one edge of $C$ in $D_H$, and, since $f$ is the only edge that crosses $e$ in $D_H$, $f$ must be in $C$, while neither $u$ nor $v$ is in $C$.

Let $D$ be a planar embedding of $G-e$ and let $B_u$ and $B_v$ be the $C$-bridges in $G-e$ containing $u$ and $v$, respectively.  If $B_u=B_v$, then there is a $C$-avoiding $uv$-path $P$ in $B_u$.  Then $P+e$ and $C$ are disjoint cycles that show $e$ and $f$ are separated by cycles in $G$, a contradiction.  

Thus, $B_u\ne B_v$.  Since $B_u^H\subseteq B_u$ and $B_v^H\subseteq B_v$, $B_u$ and $B_v$ are overlapping $C$-bridges in $G$.  In particular, they are in different sides of $C$ in $D$.

Let $G_u$ and $G_v$ be the subgraphs of $G$ embedded in the side of $C$ containing $u$ and $v$ in $D$, respectively. Our goal is to prove there exist planar embeddings of $G_u$ and $G_v$ such that, in both cases, $C$ bounds a face and, in $G_u$, $u$ and $f$ are cofacial, while in $G_v$, $v$ and $f$ are cofacial.  This obviously implies that $G$ has a 1-drawing in which $e$ and $f$ are the crossing pair.  We prove the result for $G_u$ and the result for $G_v$ follows by symmetry. 

\begin{figure}
	\centering
	\includegraphics[scale=0.7]{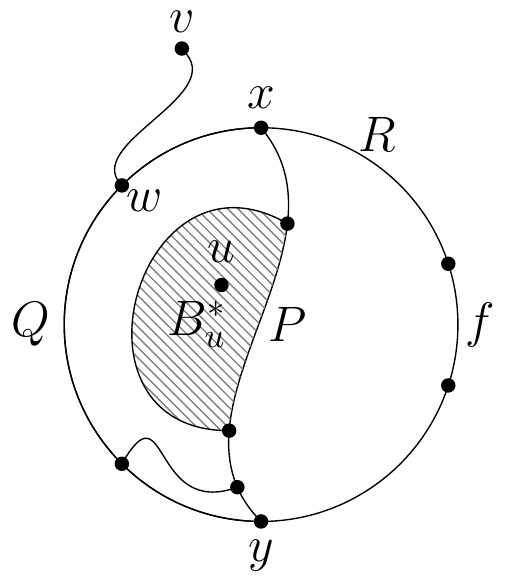}
	\caption{An abstract representation of the elements introduced for the proof of Theorem \ref{thm:1cgc}.}
	\label{fig:propria}
\end{figure}

We refer the reader to Figure \ref{fig:propria} for a visual aid in what follows. Let $F$ be the face of $D[C\cup B_u]$ incident with $f$ but not bounded by $C$.  There is a unique cycle $C_f$ in $C\cup B_u$ that contains $f$ and is contained in the boundary of $F$.  To see that $C_f$ exists, the cycle $C$ is one cycle that has $F$ on one side and contains $f$.  Choose $C_f$ to be (the) one that has a minimal subgraph of $G$ on the side containing $F$. This choice implies that every $C_f$-bridge embedded on the side containing $F$ has at most one attachment.

Let $R$ be the component of $C \cap C_f$ including $f$. As $B_v$ overlaps $B_u$ on $C$, $B_u$ has at least two attachments on $C$.  Thus, $R$ is the subpath of $C$ that contains $f$ and joins two consecutive attachments of $B_u$. Let $x$ and $y$ be the two ends of $R$, let $P$ be the other $xy$-path in $C_f$, and let $Q$ be the other $xy$-subpath of $C$.  

We note that every $(C\cup P)$-bridge in $C\cup B_u$ is a $(P\cup Q)$-bridge and, if a $(C\cup P)$-bridge $B$ in $C\cup B_u$ is contained in the face of $C\cup P$ bounded by $P\cup R$, then $Att(B)$ is a single vertex in $P$.   

To show that $G_u$ has an embedding in which $u$ and $f$ are cofacial and $C$ bounds a face, we shall show that $C\cup B_u$ has such an embedding.   If $u$ is in $P$, then $D$ is already such an embedding, so we assume $u$ is not in $P$.  Recall that, since $C$ detaches $u$ and $v$ in $G-e$, $u$ is not in $C$.  

Therefore, there is a $(C\cup P)$-bridge $B^*_u$ in $C\cup B_u$ containing $u$.  Because $x$ and $y$ are consecutive attachments of $B_u$ on $C$ and $B_v$ overlaps $B_u$ on $C$, $B_v$ has an attachment $w$ in $Q-\{x,y\}$.  Thus, 

\medskip (*) there is a $C$-avoiding $vw$-path $P_v$ in $B_v$.

\medskip
Suppose $B^*_u$ has an attachment $z$ in $Q-\{x,y\}$; let $P^*$ be a $(C\cup P)$-avoiding $uz$-path in $B^*_u$.  Letting $Q[w,z]$ denote the $wz$-subpath of $Q$,  $(P_v\cup Q[w,z]\cup P^*)+e$ and $P\cup R$ are disjoint cycles in $G$ containing $e$ and $f$, respectively.  This contradicts the assumption that $e$ and $f$ are not separated by cycles.  Thus, 

\medskip (**) $Att(B^*_u)$ is contained in $P$.

\medskip
To show that $C\cup B_u$ has an embedding in which $u$ and $f$ are cofacial and $C$ bounds a face, we first show that $C\cup P\cup  B^*_u$ has such an embedding.
For the sake of contradiction, suppose that there is no such embedding of $C\cup P\cup B^*_u$.  

By Corollary \ref{cor:tuttefc}, there is a cycle $C'$ in $C\cup P\cup B^*_u$  that detaches $u$ and $f$.  Let $B'_u$ and $B'_f$ be the $C'$-bridges containing $u$ and $f$, respectively.  Note that $B'_u\subseteq B^*_u$ and $B'_u$ and $B'_f$ are overlapping $C'$-bridges.

To see that $Q\subseteq C'$, suppose $C'\subseteq B^*_u\cup P$.  Let $z$ be any attachment of $B'_f$.  Then there is a $C'$-avoiding path $P'_z$ in $B'_f-f$ joining an end of $f$ to $z$.  Choose the labelling of $x$ and $y$ so that $x$ is in $P'_z$.  

As we saw in (*) above, there is a $C$-avoiding  $vw$-path $P_v$ in $B_v$.  Extend this path to $x$ using the $wx$-subpath of $Q$, and finally extend to $z$ using the $xz$-subpath of $P'_z$.  The result is a $C'$-avoiding $vz$-path in $G-\{e,f\}$. 

Consider the graph $K$ consisting of $(C-f)\cup B^*_u\cup B_v$.  We have shown that the $C'$-bridge $B'_v$ in $K$ that contains $B_v$ has all the attachments of $B'_f$ as attachments.  Since $B'_f$ overlaps $B'_u$, it follows that $B'_v$ overlaps $B'_u$.

Therefore, Observation \ref{obs:overlapping} shows that $K+e$ is not planar.  However, $K+e$ is contained in $G-f$, contradicting the hypothesis that $G-f$ is planar.  Thus, $Q\subseteq C'$, so $C'$ consists of $Q$ and an $xy$-path contained in $B^*_u\cup P$.

It follows that $B'_f=R$ and, therefore, $B'_u$ has an attachment in each component of $C'-\{x,y\}$; in particular, we have the contradiction  that $B'_u$, and therefore $B^*_u$, has an attachment in $Q-\{x,y\}$, contradicting (**).

Consequently, $u$ and $f$ are cofacial in an embedding $D_1$ of $C\cup P\cup B^*_u$ in which $C$ bounds a face.  We first extend this to an embedding $D_2$ of $C\cup B_u$ in which $u$ and $f$ are cofacial and $C$ bounds a face.  

Let $B$ be a $(C\cup P\cup B_u^*)$-bridge in $C\cup B_u$.   If $B$ is embedded by $D$ in the face of $C\cup P$ bounded by $P\cup R$, then $B$ has only one attachment and this is in $P$.  Thus, it may be added to $D_1$ while retaining the facts that $u$ and $f$ are cofacial and $C$ bounds a face.

The remaining $(C\cup P\cup B_u^*)$-bridges are embedded by $D$ in the face of $C\cup P$ bounded  by $P\cup Q$.  These may all be added in that same face of $D_1$, completing the required embedding $D_2$.

Finally, we  add all the remaining $(C\cup P)$-bridges in $G_u$ to maintain the facts that $u$ and $f$ are cofacial and $C$ bounds a face.  If $B$ is embedded by $D$ in the face of $C\cup P$ bounded by $P\cup Q$, then we use that embedding in $D_2$.  

\begin{figure}
	\centering
	\subcaptionbox{\label{fig:frame2}}
	[.45\linewidth]{\includegraphics[scale=0.7]{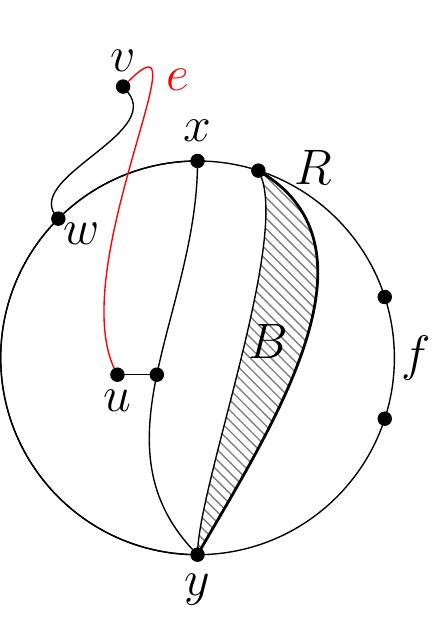}}
	\subcaptionbox{\label{fig:frame}}
	[.45\linewidth]{\includegraphics[scale=0.7]{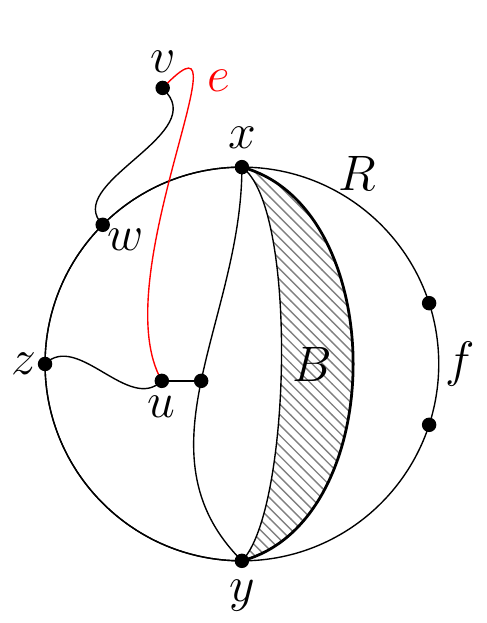}}
		\caption{}
		\label{fig:frameall}
\end{figure}

In the last case, suppose  $B$ is embedded by $D$ in the face of $C\cup P$ bounded by $P\cup R$ (see Figures \ref{fig:frame2} and \ref{fig:frame}).  If $B$ has an attachment in $P-\{x,y\}$, then $B\subseteq B_u$ is already embedded in $D_2$.  Therefore, $Att(B)$ is contained in $R$.

If all the attachments of $B$ are in one of the components of $R-f$, then we can use the embedding of $B$ in $D$ to add $B$ to $D_2$.
Thus, we may assume that $B$ has an attachment in each component of $R-f$. 

If either $B_u$ has an attachment $z$ in $Q-\{x,y\}$ or $B$ has an attachment not in $\{x,y\}$, then there are vertex-disjoint paths $P_u\subseteq B_u$ from $u$ to $Q$, and $P'$ in $B$ joining attachments  in each component of $R-f$ such that $V(P_u)\cap V(P')=\emptyset$.
Thus, the cycle in $G$ containing $e$, $P_u$, $P_v$, and the $zw$-subpath of $Q$ is disjoint from the cycle containing $f$, $P’$, and a subpath of $R$.  That is, $e$ and $f$ are separated by cycles, a contradiction.

Therefore $Att(B)=Att(B_u)=\{x,y\}$. This shows that we can add $B$ in the face of $D_2[C\cup P]$ bounded by $P\cup Q$. This completes the embedding of $G_u$ for which $u$ and $f$ are cofacial and $C$ bounds a face.

\end{proof}

\section{Final remarks}

Let $G$ be a nonplanar graph. We say that a pair of edges $\{e,f\}$ of $G$ is a \textit{potential crossing pair} if,
 for every Kuratowski subgraph $H$ of $G$, $e$ and $f$ is a crossing pair of $H$. Clearly, a crossing pair is also a potential crossing pair and, as Theorem~\ref{thm:1cgc} (ii) shows, a potential crossing pair not separated by cycles is also a crossing pair. This raises the question on the existence of potential crossing pairs on graphs with crossing number at least 2. Given Theorem~\ref{thm:1cgc}, such a pair would necessarily be separated by cycles. 
 
 \begin{conjecture}
 A non-planar graph $G$ has crossing number at least 2 if and only if it does not have a potential crossing pair.
 \end{conjecture}
 
 \v Sir\'a\v n's analysis \cite{siran} is instructive here.  Deleting any two edges of $K_6$ results in a graph with crossing number 1.  Thus, there is necessarily a potential crossing pair.  Deleting any edge $e$ of $K_6$ results in a graph with crossing number 2, but a simple case check shows there is no pair of edges that is in every $K_{3,3}$.
That is, $K_6-e$ has no potential crossing pair.

Regarding the proof of Theorem~\ref{thm:1cgc}, we are aware of a shorter proof that shows that (ii) and (i) are equivalent (see Silva's thesis~\cite[Theorem 3.12]{thesis}). This proof uses the famous Two Disjoint Paths theorem, proved independently by many authors~\cite{dpg,psutp,2lg,gmdcp}. The proof relies on the fact that a pair of edges $\{e,f\}$ of a graph $G$ is a crossing pair if and only if $G-e-f$ has a planar embedding such that the ends of $e$ and $f$ alternate in a face of the embedding. The proof of the equivalence of (i) and (ii) using the Two Disjoint Paths Theorem is shorter and somewhat straightforward. However, proving that (iii) implies (ii) seems significantly more complicated and we were unable to make it work using only the Two Disjoint Paths Theorem.

\section*{Acknowledgements}
The first author was supported by FAPESP Proc. 2015/04385-0, 2014/14375-9 and 2015/11937-9, CNPq Proc. 311373/2015-1. The second author was supported by CONACyT and by the ISTplus Fellowship. The third author was supported by NSERC Grant No.\ 41705-2014 057082. The fourth author was supported by CNPq Proc. 311373/2015-1, CNPq Proc. 425340/2016-3 and FAPESP Proc. 2015/11937-9. This project has received funding from the European Union's Horizon 2020 research and innovation
programme under the Marie Skłodowska-Curie grant agreement No 754411.

 \printbibliography

\end{document}